\documentclass[a4paper,12pt,twoside]{article}

\input{epsf.tex}
\usepackage{amsthm}
\usepackage{amsmath}
\usepackage{amssymb}
\usepackage{pstricks,pst-node,multido}
\usepackage{ragged2e}

\newcommand{\mytitle}{Generalisations of the Tits representation}

\DeclareMathOperator{\sym}{sym}
\DeclareMathOperator{\supp}{supp}
\DeclareMathOperator{\GL}{GL}

\newcommand{\td}{\mathbin{\triangledown}}
\newcommand{\gsa}{gsa}

\psset{dotsize=4pt}

\newcommand{\dashed}{\psset{linestyle=dashed, dash=3pt 2pt}}

\newcommand\sur{\mathrel{\to\kern-1.8ex\to}}

\newcommand\ra{\rightarrow}

\newcommand{\<}{\langle} 
\renewcommand{\>}{\rangle}

\newcommand{\be}{\begin{equation}}
\newcommand{\ee}{\end{equation}}
\newcommand{\BOX}{{}\raisebox{-.05ex}{\makebox[1em][r]{$\Box$}}}
\newcommand{\block}{\hfill\BOX\par}
\renewcommand{\qed}{\block}
\newcommand{\col}{\text{\upshape :\ }}

\usepackage[left=35mm, right=35mm, top=35mm, bottom=15mm, a4paper, includefoot]{geometry}

\newcommand{\headsize}{\small}

\usepackage{fancyhdr}          
\makeatletter

\renewcommand\tableofcontents{%
    \subsection*{\contentsname
          \@mkboth{%
          \MakeUppercase\contentsname}{\MakeUppercase\contentsname}}%
    \@starttoc{toc}%
    }

\def \@floatboxreset {%
        \reset@font
        \normalsize\sf
        \@setminipage
}

\makeatother

\makeatletter

\newcommand{\figurefont}{\sf\bfseries}
\newcommand{\mycaption}[2]
{\small{\figurefont #1.\ }\small{\figurefont\sf #2}}

\long\def\@makecaption#1#2{%
  \vskip\abovecaptionskip
  \sbox\@tempboxa{{\ \ \mycaption{#1}{#2} }}%
  \ifdim \wd\@tempboxa >\hsize
    {{\mycaption{#1}{#2}}}\par
  \else
    \global \@minipagefalse
    \hb@xt@\hsize{\hfil\box\@tempboxa\hfil}%
  \fi
  \vskip\belowcaptionskip}

\makeatother

\newcommand{\daandot} {\raisebox{.2ex}{\footnotesize$\circ$}}

\newcommand{\vardot}{\daandot}

\newcommand{\daanlabel} {\makebox[0em][r] {\vardot\hspace{1ex}}}
\newcommand{\itemn}{\item[\daanlabel]}

\newcounter{daana}

\newlength{\daanleftmargin}
\newlength{\daanrightmargin}
\newenvironment{mathlist}{\par%
\begin{list}{\makebox[0em][l]{%
\daanlabel%
\hspace*{-\daanleftmargin}%
\makebox[\columnwidth][r]
  {\refstepcounter{equation}%
  \upshape(\theequation)}}}
{\usecounter{daana}
\setcounter{daana}{\theequation}
\setlength{\rightmargin}{\daanrightmargin}
\setlength{\leftmargin}{\daanleftmargin}
\setlength{\labelwidth}{\daanleftmargin}
\setlength{\labelsep}{0ex}
\setlength{\itemsep}{.8ex}
\setlength{\topsep}{.8ex}
\setlength{\parsep}{0ex}
\setlength{\listparindent}{\parindent}
}}
{\end{list}}

\newcommand{\bor}[1]
{\hfill\makebox[0mm][l]{\makebox[\daanrightmargin][r]{#1}}}

\newenvironment{mathlistn}
{\begin{list}{\daandot}{
\setlength{\itemsep}{0ex}
\setlength{\topsep}{.8ex}
\setlength{\leftmargin}{\daanleftmargin}
\setlength{\labelsep}{.8ex}
\setlength{\labelwidth}{2ex}}}
{\end{list}}

\newcommand{\wild}{}

\newcommand{\thspace}{4mm}
\newcommand{\afterblock}{\vspace{\thspace}}

\newtheoremstyle{slanted}
  {\thspace}
  {\thspace}
  {\sl}
  {}
  {\bfseries}
  {.}
  {.5em}
  {}

\theoremstyle{slanted}

\newtheorem{prop}[equation]{Proposition}

\newtheorem{theo}[equation]{Theorem}
\newtheorem{lemm}[equation]{Lemma}
\newtheorem{coro}[equation]{Corollary}
\newtheorem{wildsl}[equation]{\wild}

\newtheoremstyle{roman}
  {\thspace}
  {\thspace}
  {\rm}
  {}
  {\it}
  {.}
  {.5em}
  {}

\theoremstyle{roman}

\newtheorem{defi}[equation]{Definition}
\newtheorem{rema}[equation]{Remark}
\newtheorem{ques}[equation]{Question}
\newtheorem{exam}[equation]{Example}

\newtheorem{wildrm}[equation]{\wild}

\newcommand{\daanmath}{\mathbb}
\newcommand{\zz}{\daanmath{Z}}

\newcommand{\rr}{\daanmath{R}}

\begin{document} 
%

{\Huge
\begin{center} \mytitle \\
{\normalsize Daan Krammer} \\[-1ex]
{\normalsize 9 August 2007} \\[3ex]
\end{center}}

\thispagestyle{empty}

\begin{abstract} We construct a group $K_n$ with properties similar to infinite Coxeter groups. In particular, it has a geometric representation featuring hyperplanes and simplicial chambers. The generators of $K_n$ are given by $2$-element subsets of $\{0,\ldots,n\}$. We give some easy combinatorial results on the finite residues of $K_n$.

Math classification 2000: primary 52B30, secondary 20F55, 22E40.
\end{abstract}


\section{Introduction}

A {\em Coxeter group\,} is a group $W$ presented with generating set $S$ and relations $s^2$ for all $s\in S$ and at most one relation $(st)^{m(s,t)}$  for every pair $\{s,t\}\subset S$ (where $m(s,t)=m(t,s)$). It is known that then the natural map $S\ra W$ is injective; we think of it as an inclusion. We call the pair $(W,S)$ a {\em Coxeter system}.

We generalise this as follows. For any set $S$, let $F_S$ denote the free monoid on $S$. A {\em fully coloured graph\,} is a triple $(V,S,m)$ where $V,S$ are sets, $m\col V\times S\times S\ra \zz_{\geq1}\cup\{\infty\}$ is a map, and an {\em action\,} $V\times F_S\ra V$ written $(v,g)\mapsto vg$ is specified, satisfying the following.
\begin{mathlistn}

\item For all $v\in V$, $s\in S$ we have $(vs)s=v$. 

\item Let $v\in V$, $s,t\in S$. Then
$m(v;s,t)=1$ if and only if $s=t$. Moreover $m(v;s,t)=m(v;t,s)$ and $m(v;s,t)=m(vs;s,t)$. Also, if $k:=m(v;s,t)$ is finite then $v(st)^k=v$.

\item The set $V$ is universal. That is, let $(V',S,m')$ satisfy the above too and let $f\col V'\ra V$ be a map satisfying (i) $(fv)s=f(vs)$ for all $v,s$; (ii) $m'(v;s,t)=m(fv;s,t)$ for all $v,s,t$. Then the restriction of $f$ to every $F_S$-orbit in $V'$ is injective.

\end{mathlistn}

Every Coxeter system $(W,S)$ gives rise to a {\em Coxeter fully coloured graph\,} $(W,S,m)$ where one defines $m(w;s,t)$ to be the order of $st$ and the action $W\times S\ra W$ to be multiplication.

Let $(V,S,m)$ be a fully coloured graph. For $I\subset S$, an {\em $I$-residue\,} is a subset of $V$ of the form $\{vg\mid g\in F_I\}$. Contrary to the Coxeter case, it may happen that an $\{s,t\}$-residue $R\subset V$ containing $v$ is such that $\#R\neq 2m(v;s,t)$; see remark~\ref{tr62}(b) for an example of this.

A {\em generalised simplicial hyperplane arrangement\,} (\gsa) consists of an open convex cone $U^0$ in a finite dimensional real vector space $Q$ together with a set $\cal A$ of hyperplanes%
\footnote{A hyperplane is a $1$-codimensional linear subspace.}
which is locally finite%
\footnote{Locally finite in $U^0$ means that every compact subset of $U^0$ meets finitely hyperplanes from $\cal A$.}
in $U^0$ and such that every (closed) chamber%
\footnote{A chamber is the closure of a connected component of $U^0\smallsetminus(\cup {\cal A})$.}
is simplicial%
\footnote{A simplicial chamber is one of the form $\{x\in Q\mid f_i(x)\geq 0$ for all $i\}$ for some independent set of linear maps $\{f_i\mid i\in I\}$.}
and such that for every chamber $C$ and every boundary point $c\in\partial C$ there is a chamber $D\neq C$ containing $c$.
For two chambers $C,D$, let $d(C,D)$ be the number of hyperplanes $H\in\cal A$ separating $C^0$ from $D^0$ ($0$ is relative interior). A {\em cell\,} is an intersection of chambers (we use a different but equivalent definition in the main text). If $U^0=Q$ then $\cal A$ is finite and we obtain a simplicial hyperplane arrangement.

Every \gsa\ $\cal A$ gives rise to a fully coloured graph $G({\cal A})=(V,S,m)$ which is defined by the following. We let $V$ be the set of chambers. Let $n$ denote the codimension of the smallest cell. There is a unique equivalence relation on the set of $(n-1)$-codimensional cells of $n$ equivalence classes such that no two equivalent $(n-1)$-codimensional cells are in one chamber. Let $S$ be the set of equivalence classes. Let the action $V\times S\ra V$ be such that, for all $C,s$, one has $d(C,Cs)=1$ and no $(n-1)$-codimensional cell in $s$ is contained in $C\cap Cs$. Let $m(C;s,t)$ be half the number of chambers containing $C\cap Cs\cap Ct$.

A {\em realisation\,} of a coloured graph $\Gamma$ is a \gsa\ $\cal A$ such that $G({\cal A})\cong \Gamma$. For a fully coloured graph, a realisation may not exist, and it may not be unique up to isomorphism if it exists.

A celebrated result by Tits \cite[section~5.4.4]{bou}, \cite{vin}, \cite[section~5.13]{hum} states that every Coxeter group (seen as a fully coloured graph) can be realised. Moreover, the realisation can be chosen to be covariant under some $W$-action, so that we have a faithful linear representation $W\ra \GL(Q)$.

It is not hard to generalise Tits's result, with essentially the same proof, to theorem~\ref{6} which roughly states that something is a \gsa\ as soon as its $3$-residues are, or even if certain small parts of $3$-residues are. The challenge lies in finding interesting examples where theorem~\ref{6} can be used to prove that something is a \gsa. To arrive at such examples, one has to overcome two obstacles which are easy in the case of Coxeter groups: (a) to find a fully coloured graph; and (b) to find a realisation of it.

We give a modest partial solution to part (b) as follows. Define a $(2,3)$-graph to be a fully coloured graph $(V,S,m)$ such that $m(v;s,t)\in \{2,3\}$ for all $v,s,t$ and such that all $3$-residues are either of the common Coxeter type, or of type $A(3,7)$ defined in figure~\ref{tr61}(b). (In the main text we use a different but equivalent definition of $(2,3)$-graphs.) Then every $(2,3)$-graph is realisable. Without much more effort one proves a more general result involving so-called $(2,3,\infty)$-graphs which we also include (see theorem~\ref{amcg6} and proposition~\ref{tr49}).

The fully coloured graph associated with a Coxeter system $(W,S)$ is a $(2,3)$-graph if and only if the Coxeter system is simply laced, that is, the order of $st$ is in $\{1,2,3\}$ for all $s,t\in S$.

For $n\geq 0$, let $K_n$ be the group presented by a set $T_n\subset K_n$ of $\binom{n+1}{2}$ generators written
\[ T_n=\Big\{ t(a,b)=\binom{a}{b} \ \Big|\ a,b\in\{0,1,\ldots,n\},\ a< b \Big\} \]
and relations $s^2$ for all $s\in T_n$ and
\[ \binom{a}{b}\binom{c}{d}\binom{a}{b}\binom{c}{d} \]
whenever $0\leq a< b\leq c< d\leq n$;
\[ \binom{a}{b}\binom{a+x}{b-y}\binom{a}{b}\binom{a+y}{b-x} \]
whenever $x,y\geq0$ and $0\leq a<a+x+y<b\leq n$; and
\[ \binom{a}{b-z} \binom{a+y}{b} \binom{a}{b-x} \binom{a+z}{b} \binom{a}{b-y} \binom{a+x}{b} \]
whenever $x,y,z> 0$ and $0\leq a\leq a+x+y+z=b\leq n$.

This construction is motivated by the observation that there exists a $K_n$-action on $\{1,\ldots,n\}$ given by
\[ \binom{a}{b}(x)=\left\{ \begin{array}{@{}lll@{}} a+b+1-x &\hspace{2mm}& \text{if $a+1\leq x\leq b$,} \\[.5ex] x && \text{otherwise.} \end{array} \right. \]

Our main result (theorem~\ref{tr55}) states that there exists a $K_n$-action on a $(2,3)$-graph $\Gamma_n=(V,S,m)$ such that the action on $V$ is simply transitive, and such that there exists a vertex $v\in V$ such that, for all $g\in K_n$, we have $d(v,gv)=1$ if and only if $g\in T_n$. We give a case-by-case proof of the theorem by looking at every $3$-residue separately.

As we remarked above, $(2,3)$-graphs are realisable. In particular, $K_n$ has a geometric representation much as Coxeter groups have.

Residues of realisable fully coloured graphs (seen as fully coloured graphs themselves) are again realisable. A residue in a Coxeter fully coloured graph is again Coxeter. Contrary to this, a residue in $\Gamma_n$ is not necessarily isomorphic to any $\Gamma_k$. We call such residues {\em admissible graphs\,} and we study them on a par with $\Gamma_n$ itself.

Among the $(2,3)$-graphs the finite ones seem most interesting. We list the irreducible rank~$4$ $(2,3)$-graphs without proof in proposition~\ref{tr59}. Two of them are Coxeter and two of them are not. Both of the non-Coxeter ones are admissible. This suggests that $\Gamma_n$ may be a good source for finite $(2,3)$-graphs.

Section~\ref{tr56} introduces coloured graphs and 
proves Tits's result in our wider setting. This also generalises the fully coloured graphs mentioned above. In section~\ref{tr57} we define $(2,3,\infty)$-graphs (which in the main text are by definition realisable) and classify them locally. In section~\ref{tr58} we study the group $K_n$ and its relation with $(2,3)$-graphs.

\section{Realisations of 
coloured graphs} \label{tr56}

A {\em partial map\,} $f\col A\ra B$ (of sets, say) consists of a subset $D\subset A$ and a map $D\ra B$. We call $D$ the domain of $f$ and we say that $f(a)$ is not defined unless $a\in D$. A statement such as ``$f(a)$ is positive'' implies in particular that $f(a)$ is defined.

For a set $S$, let $F_S$ be the free monoid on $S$. We consider $S$ to be a subset of $F_S$. If $S\subset T$ then $F_S\subset F_T$.

\begin{defi} A {\em coloured graph\,} is a tuple $(V,S,A,m)$ with the following properties. Firstly, $V$ is a set (of vertices) and $S$ is a set (of colours). We have a partial action $A\col V\times F_S\ra V$ written $(v,g)\mapsto vg$, that is, a partial map such that for all $v\in V$, $g,h\in F_S$, if $v(gh)$ or $(vg)h$ is defined then so is the other, and they are equal. Usually we omit $A$ from the notation. We have a partial map $m\col V\times S\times S\ra \zz_{\geq1}\cup{\infty}$ such that $m(v;s,t)$ is defined if and only if, for all $k>0$, the vertices $v(st)^k$ and $v(ts)^k$ are defined. Moreover, the following hold.
\begin{mathlist}

\itemn If $vs$ is defined ($v\in V$, $s\in S$) then $(vs)s=v$. 

\itemn We have $m(v;s,t)=1$ if and only if $s=t$ and $vs$ is defined.  (Recall that if $vs$ is not defined then neither is $m(v;s,t)$.)

\item \label{tr53} Suppose that $m(v;s,t)$ is defined. Then $m(v;s,t)=m(v;t,s)$ and $m(v;s,t)=m(vs;s,t)$.

\item \label{tr23} If $k:=m(v;s,t)$ is defined and finite then $v(st)^k=v$.

\item \label{tr24} The set $V$ is universal. That is, let $(V',S,m')$ satisfy the above too and let $f\col V'\ra V$ be a map satisfying (i) $(fv)s=f(vs)$ for all $(v,s)\in V'\times S$ such that at least one side is defined; (ii) whenever $m'(v;s,t)$ is defined, it equals  $m(fv;s,t)$. Then the restriction of $f$ to every $F_S$-orbit in $V'$ is injective.

\end{mathlist}A {\em fully coloured graph\,} is a coloured graph $(V,S,m)$ such that the action $V\times S\ra V$ (hence $m$) is everywhere defined. A direct definition of fully coloured graphs was given in the introduction.\qed\end{defi}

Let $(V,S,m)$ be a coloured graph and let $I\subset S$. An {\em $I$-residue\,} is a subset of $V$ of the form $\{vg\mid g\in F_I\}$ where $v\in V$ (of course, only the defined vertices are included). We call it also an {\em $r$-residue\,} if $r=\#I$.

Let $R$ be the $\{s,t\}$-residue through $v$. It follows from (\ref{tr53}) that $m(v;s,t)$ depends only on $(R,s,t)$ (if it is defined). We write it as $m(R;s,t)$ accordingly.

As explained in the introduction, every Coxeter system gives rise to a fully coloured graph.

\begin{rema} \label{tr62} Let $(V,S,m)$ be a coloured graph. 

(a). There is an equivalence relation on $V$ with two equivalence classes such that $v$, $vs$ are not equivalent for all $v\in V$, $s\in S$. In particular, $v\neq vs$. This follows from the universality (\ref{tr24}) and the fact that the relations (\ref{tr23}) have even length.

(b). Let $R$ be an $\{s,t\}$-residue. Then $\#R$ divides $2m(R;s,t)$, but it may happen that $\#R\neq 2m(R;s,t)$, as is shown by the following example. Put $V=(\zz/2)^3$ and $S=\{r,s,t\}\subset V$ where $r=(1,0,0)$, $s=(0,1,0)$, $t=(0,0,1)$. Let $S$ act on $V$ by right multiplication. Define $m(v;a,b)=2$ for all $v,a,b$ except if $\{a,b\}=\{s,t\}$ and $v\in R:=\<s,t\>$ in which case we put $m(v;s,t)=4$. Then $(V,S,m)$ is a coloured graph but $2m(R;s,t)=8\neq4=\#R$.

(c). From section~\ref{tr57} we shall only consider {\em fully\,} coloured graphs. We take the opportunity to prove the results of this section in the slightly more general setting of partial maps.
\end{rema}

Let $Q$ be a real vector space. A {\em hyperplane\,} in $Q$ is a $1$-codimensional linear subspace. An {\em open (respectively, closed) half-space\,} is a subset of $Q$ of the form $f^{-1}(\rr_{>0})$ (respectively, $f^{-1}(\rr_{\geq 0})$) where $f\col Q\ra\rr$ is a nonzero linear map. If $H$ is one of the above half-spaces, then the {\em boundary\,} $\partial H$ is defined to be $f^{-1}(0)$.

We call a coloured graph $(V,S,m)$ {\em connected\,} if for all $v,w\in V$, there exists $g\in F_S$ such that $vg=w$.

\begin{defi} Let $\Gamma=(V,S,m)$ be a connected coloured graph. A {\em realisation\,} of $\Gamma$ consists of the data (\ref{tr4})--(\ref{tr6}) satisfying properties (\ref{tr1})--(\ref{tr3}) below.
\begin{mathlist}

\item \label{tr4} For every $v\in V$ a real vector space $P(v)$ with basis $\{p(v,s)\mid s\in S\}$ (a set in bijection with $S$).

\item \label{tr6} Whenever $w:=vs$ is defined ($v\in V$, $s\in S$) an isomorphism 
\[ \phi_{v,s}\col P(v)\ra P(w) \] 
such that $p(v,t)\,\phi_{vs}=p(w,t)$ for all $t\in S\smallsetminus{s}$.

\item \label{tr1} Let $Q$ denote the quotient of the disjoint union $\sqcup_{v\in V}P(v)$ by the smallest equivalence relation $\equiv$ such that $x\, \phi_{vs}\equiv x$ for all $v,s$ and all $x\in P(v)$. Then the natural map $P(v)\ra Q$ is bijective for one hence all $v\in V$.
\end{mathlist}
Note that the condition (\ref{tr1}) is equivalent to $\phi_{v_1s_1}\cdots \phi_{v_ns_n}=1$ (indices in $\zz/n$) whenever $v_is_i=v_{i+1}$ for all $i$. It is sufficient for this to hold for $\#\{s_1,\ldots,s_n\}=2$, by (\ref{tr24}).

The image in $Q$ of $p(v,s)$ is written $q(v,s)$. It follows from (\ref{tr1}) that $Q$ is a real vector space with basis $\{q(v,s)\mid s\in S\}$ (a set in bijection with $S$) whenever $v\in V$. For $v\in V$ we define the {\em chamber\,} $C(v)=\sum_{s\in S}\rr_{\geq 0} \, q(v,s)$. 

\begin{mathlist}

\item \label{tr5} We have 
$C(v)^0\cap C(vs)^0=\varnothing$ for all $v\in V$, $s\in S$,  where $0$ denotes the relative interior.

\item \label{tr3} Let $R\subset V$ be an $\{s,t\}$-residue, $s\neq t$. Suppose that $X=\cap_{v\in R}C(v)$ has codimension $2$, that is, $\#R\geq2$.

If $k=m(R;s,t)$ is defined and finite then there exist $k$ (distinct) hyperplanes in $Q$ containing $X$ such that every component of the complement of these hyperplanes meets $C(v)$ for a unique $v\in R$. In particular, $\#R=2m(R;s,t)$.

If $m(R;s,t)$ is infinite or not defined then $\cup_{v\in R}C(v)$ is contained in some closed half-space whose boundary contains~$X$.\bor{\BOX}

\end{mathlist}
\end{defi}

Suppose $vs=w$ ($v\in V$, $s\in S$). Then there are unique $c_t\in\rr$ ($t\in S$) such that
\[ q(w,s)= \sum_{t\in S} c_t\, q(v,t). \]
Now (\ref{tr5}) is equivalent to $c_s<0$.

\begin{exam} It is not hard to show that every Coxeter coloured graph admits a (covariant) realisation determined by 
\be p(v,s)\, \phi_{v,s}= -p(vs,s)+\sum_{t\in S\smallsetminus\{s\}} 2\, \cos\frac{\pi}{m(s,t)}\, p(vs,t). \label{tr65} \ee
See \cite[section~5.4.3]{bou}, \cite[section~5.3]{hum}, \cite{vin}.\end{exam}

\begin{rema} Suppose that the coloured graph $(V,S,m)$ admits a realisation. Let $v\in V$ and let $s,t\in S$ be distinct. If $m(v;s,t)$ is defined (but possibly infinite) then the $\{s,t\}$-residue through $v$ has $2m(v;s,t)$ elements. This follows immediately from (\ref{tr3}). In particular, $vs\neq vt$.

In the case of Coxeter groups, this is the usual proof that the order of $st$ equals $m(s,t)$ rather than a proper divisor of it.\end{rema}

Let $(V,S,m)$ be a connected coloured graph. For $v,w\in V$, define $d(v,w)$ to be the least $k\geq0$ such that there are $s_1,\ldots,s_k\in S$ with $vs_1\cdots s_k=w$. Then $d$ is a metric. By a {\em semi-geodesic\,} we mean a tuple $(v_1,\ldots, v_n)$ of vertices such that $d(v_1,v_n)=\sum_{i=1}^{n-1} d(v_i,v_{i+1})$.

For $v\in V$, $s\in S$, we define
\[ H(v,s):=\Big\{\sum_{t\in S} c_t\, q(v,t)\ \Big|\ c_t\in\rr\text{ for all }t\in S \text{ and } c_s\geq0\Big\} \subset Q. \]
If $vs$ exists then $H(v,s)$ is the closed half-space in $Q$ containing $C(v)$ whose boundary contains $C(v)\cap C(vs)$.

In the remainder of this section, we consider a connected coloured graph with a realisation, and use the above notation. 

The proof of the following is similar to \cite[section~5.4.4]{bou}.

\begin{theo} \label{6} Let $v',w\in V$ be distinct. Let $s\in S$ and suppose that either (a) $v's$ is not defined, or (b) $v:=v's$ is defined and $(v,v',w)$ is a semi-geodesic. Then $C(w)\subset H(v',s)$.
\end{theo}

\begin{proof} Induction on $n=d(v',w)$. For $n=0$ it is trivial. If $n\geq1$, let $v''=v't$ ($t\in S$) be a neighbour of $v'$ such that $(v',v'',w)$ is a semi-geodesic. Note that $s\neq t$. In case (b), we have $v''\neq v$.

Let $R$ be the $\{s,t\}$-residue through $v'$. For $a,b\in R$, let $d_0(a,b)$ be the least $k\geq 0$ such that there exist $s_1,\ldots,s_k\in\{s,t\}$ with $b=as_1\cdots s_k$. So $d_0(a,b)\geq d(a,b)$.

Let $A$ denote the set of those $a\in R$ for which $d(v',w)=d_0(v',a)+d(a,w)$. Let $x\in A$ be an element with $d(x,w)$ minimal.

We have $\# R\geq2$ because $v',v''\in R$. Let $y\in R$ be a neighbour of $x$, that is, $d_0(x,y)=1$.

We claim that $(y,x,w)$ is a semi-geodesic. If not, we would have $d(w,y)=d(w,x)-1$ and hence
\begin{align*} d(w,v') & \leq d(w,y) + d(y,v') \leq d(w,y) + d_0(y,v') \\
 & = (d(w,x)-1) + d_0(y,v') \\ & \leq d(w,x)-1+d_0(x,v')+1 = d(w,v').  \end{align*}
So equality holds throughout, forcing $d(w,v') = d(w,y) + d_0(y,v')$, and therefore $y\in A$, contrary to $d(w,y)<d(w,x)$. 

Note that $v''\in A$, whence $d(w,x)\leq d(w,v'')<d(w,v')$. Therefore we may apply the induction hypothesis to the triples $(x,w,r)$ for $r\in\{s,t\}$. We find that
\be \label{tr41} C(w)\subset H(x,s)\cap H(x,t). \ee
In case (a) we have $H(x,s)\cap H(x,t)\subset H(v',s)$ so by (\ref{tr41}) we find $C(w)\subset H(v',s)$ as required. Suppose now that we're in case (b). Then $d_0(x,v)>d_0(x,v')$, since otherwise
\begin{align*}  d(w,v) &\leq d(w,x) + d(x,v) \leq d(w,x)+d_0(x,v) \\
& < d(w,x) + d_0(x,v') = d(w,v'),  \end{align*}
a contradiction. By (\ref{tr24}), this shows that $H(x,s)\cap H(x,t)\subset H(v',s)$ and, on combining with (\ref{tr41}) as before, $C(w)\subset H(v',s)$.\qed\end{proof}

\begin{coro} \label{tr42} If $v,w\in V$ are distinct then $C(v)^0\cap C(w)^0=\varnothing$.
\end{coro}

\begin{proof} Let $(v,v',w)$ be a semi-geodesic with $v'=vs$, $s\in S$. Apply theorem~\ref{6}.\qed\end{proof}

\afterblock A {\em cell\,} is a set of the form $\sum_{s\in I} \rr_{\geq 0}\, q(v,s)$ (which is $\{0\}$ if $I=\varnothing$) for $v\in V$, $I\subset S$.

\begin{coro} \label{7} Let $X,Y$ be distinct cells. Then $X^0\cap Y^0=\emptyset$.
\end{coro}

\begin{proof} Let $v,w$ be vertices such that $X\subset C(v)$, $Y\subset C(w)$ with $n=d(v,w)$ minimal. (We don't assume that $X$ is a ``face'' of $C(v)$ or $Y$ is of $C(w)$.) If $n=0$ it is trivial so suppose $n>0$. Let $v'=vs$ be a neighbour of $v$ such that $(v,v',w)$ is a semi-geodesic. Then $X\not\subset C(v')$ by minimality of $n$. So $X^0\cap H(v',s)=\varnothing$. We also have $Y\subset C(w)\subset H(v',s)$ so $X^0\cap Y^0=\varnothing$.\qed\end{proof}

\afterblock
The union of all $C(v)$ is denoted $U$.
 
\begin{coro} \label{9} The following hold.

{\upshape (a).} $U$ is convex.

{\upshape (b).} For all $x,y\in U$, the line segment $[x,y]:=\{tx+(1-t)y\mid 0\leq t\leq 1\}$ meets finitely many cells of $U$.
\end{coro}

\begin{proof} By corollary~\ref{7} we can prove parts (a) and (b) at once by showing that for all $x,y\in U$, the line segment $[x,y]$ is contained in the union of finitely many cells. Let $v,w$ be vertices with $x\in C(v)$, $y\in C(w)$, $n=d(v,w)$ minimal. Induction on $n$. If $n=0$ it is trivial. If $n>0$, write $[x,y]\cap C(v)=[x,z]$. Since $y\not\in C(v)$, we have $y\not\in H(v,s)$ for some $s\in S$ with $z\in\partial H(v,s)$. Since $y\in C(w)\backslash H(v,s)$, it follows from theorem~\ref{6} that (i) $v':=vs$ is defined, and (ii) $d(v',w)<d(v,w)$. Since $z\in C(v')$, the segment $[z,y]$ is contained in finitely many cells by induction. Moreover, $[x,z]$ is clearly contained in finitely many cells. This proves the induction step which finishes the proof.\qed\end{proof}

\section{$(2,3,\infty)$-Graphs} \label{tr57}

From now on, all our coloured graphs are fully coloured.

\begin{defi} \label{tr22} A {\em $(2,3,\infty)$-graph\,} is a connected fully coloured graph $(V,S,m)$ which admits a (necessarily essentially unique) realisation (\ref{tr4})--(\ref{tr3}) with the following properties.
\begin{mathlist} 

\item \label{tr21} We have $m(v;s,t)\in\{2,3,\infty\}$ for all $v,s,t$.

\item \label{tr18} {\RaggedRight We define a bijection $N\col\{2,3,\infty\}\ra\{0,1,2\}$ by $N(2)=0$, $N(3)=1$, $N(\infty)=2$. Equivalently, $N(k)=2\cos(\pi/k)$. We put $n(v;s,t):=N\big(m(v;s,t)\big)$ and $n(R;s,t)=n(v;s,t)$ if $R$ is the $\{s,t\}$-residue through $v$.

}Suppose $vs=w$ ($v\in V$, $s\in S$). Then
\begin{align*} p(v,s)\, \phi_{vs} &=-p(w,s)+\sum_{t\in S\smallsetminus\{s\}}n(v;s,t)\, p(w,t) \\*
&=-p(w,s)+\sum_{t\in S\smallsetminus\{s\}}\cos\frac{\pi}{m(v;s,t)} \, p(w,t) \end{align*}
Compare with (\ref{tr65}).
\end{mathlist}
The realisation with these properties is called the {\em standard realisation\,} in order to distinguish it from other realisations, if any. Note that the uniqueness of the standard realisation follows immediately from (\ref{tr18}).\qed\end{defi}

Our next aim is to provide an explicit criterion for a coloured graph to be a $(2,3,\infty)$-graph. We need the notion of structure sequence, which we shall now define (see figure~\ref{tr19}).

\begin{defi} Let $(V,S,m)$ be a coloured graph satisfying  (\ref{tr21}). Let $v\in V$, let $s,t\in S$ be distinct, and suppose that $k:=m(v;s,t)$ is defined and finite. Define $v_i$ ($i\in \zz/2k$) by $v_0=v$, $v_{2i-1}t=v_{2i}=v_{2i+1}s$ for all $i$ (see figure~\ref{tr19}). The map
\begin{align*}
\zz/2k &\longrightarrow \{0,1,2\} \\
2i &\longmapsto n(v_{2i};r,s) \\
2i+1 &\longmapsto n(v_{2i+1};r,t)
\end{align*}
is called the {\em structure sequence\,} of the $\{s,t\}$-residue $R$ through $v$. We always consider two structure sequences to be equal if they differ only by a cyclic permutation or reversal. Therefore the structure sequence is determined by $(R,s,t)$.

If $m(R;s,t)$ is infinite or undefined then we don't consider an associated structure sequence.\qed\end{defi}

\newcommand{\trcc}[1]{\pscirclebox*
[framesep=1pt]{#1}}                   

\newcommand{\trau}{%
\psset{showpoints=true, yunit=.8pt, xunit=.81pt} \SpecialCoor
\psframebox[linestyle=none, framesep=0pt]{\pspicture[.5](-116,-100)(116,100)$
\pnode(-116, 30){v5} 
\pnode( -50, 30){v6} 
\pnode(  50, 30){v7} 
\pnode( 116, 30){v8}
\pnode(-116,-30){v9} 
\pnode( -50,-30){v10} 
\pnode(  50,-30){v11} 
\pnode( 116,-30){v12}
\pnode( -83, 50){v2} 
\pnode(   0, 60){v3} 
\pnode(  83, 50){v4}
\pnode( -83,-50){v13} 
\pnode(   0,-60){v14} 
\pnode(  83,-50){v15}
\pnode(0, 100){v1} 
\pnode(0,-100){v16}
\pspolygon(v1)(v2)(v5)(v9)(v13)(v16)(v15)(v12)(v8)(v4)
\pspolygon(v3)(v6)(v10)(v14)(v11)(v7)
{\footnotesize \psset{framesep=2pt}
\ncline{v3}{v6}        \ncput{\trcc{s}}
\ncline{v14}{v10}      \ncput{\trcc{s}}
\ncline{v11}{v7}       \ncput{\trcc{s}}
\ncline{v6}{v10}       \ncput{\trcc{t}}
\ncline{v11}{v14}      \ncput{\trcc{t}}
\ncline{v7}{v3}        \ncput{\trcc{t}}
\ncline{v1}{v3}        \ncput{\trcc{r}}
\ncline{v2}{v6}        \ncput{\trcc{r}}
\ncline{v13}{v10}      \ncput{\trcc{r}}
\ncline{v16}{v14}      \ncput{\trcc{r}}
\ncline{v15}{v11}      \ncput{\trcc{r}}
\ncline{v4}{v7}        \ncput{\trcc{r}}}
\psdots(v2)(v9)(v16)(v12)(v4)
\psdots(v3)(v10)(v11)
\psset{dotstyle=o}
\psdots(v1)(v5)(v13)(v15)(v8)
\psdots(v6)(v7)(v14)
\rput{0}(0,0){R}
\rput{0}(-83,0){1}
\rput{0}( 83,0){1}
\rput{0}( 33.33,  60){0}
\rput{0}( 33.33, -60){0}
\rput{0}(-33.33, -60){0}
\rput{0}(-33.33,  60){0}
\uput[-90]{0}(v3){v_0}
\uput[-30]{0}(v6){v_1}
\uput[30]{0}(v10){v_2}
\uput[90]{0}(v14){v_3}
\uput[150]{0}(v11){v_4}
\uput[-150]{0}(v7){v_5}
$\endpspicture}}

\begin{figure}[h] \caption{Structure sequences.\label{tr19}} \vspace{-11mm}
\begin{flalign*} \parbox{.47\textwidth}{\RaggedRight This picture shows part of a $3$-residue $T$ containing an $\{s,t\}$-residue $R=\{v_i\mid i\}$ with $m(v_0;s,t)=3$. In the middle of every $2$-residue $R_i$ in $T$ meeting $R$ in an edge $\{v_i,v_{i+1}\}$ the picture shows the value of $n(R_i;s,t)$. The structure sequence for $R$ is $(0,0,1,0,0,1)$.} 
&& 
\trau
\end{flalign*}
\end{figure}

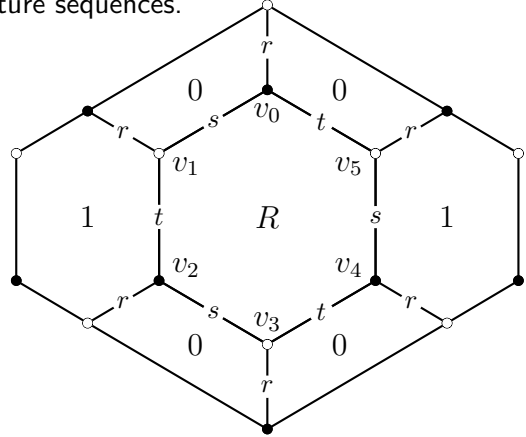
\begin{theo} \label{amcg6} {\upshape (a).} Let $\Gamma$ be a connected fully coloured graph satisfying (\ref{tr21}). Then $\Gamma$ is a $(2,3,\infty)$-graph if and only if the following hold.
\begin{mathlist}
\item \label{tr9} All structure sequences of length 4 are of the form $(n_1,n_2,n_1,n_2)$, $n_1,n_2\in\{0,1,2\}$.
\item \label{tr10} All structure sequences of length 6 are of the form $(n_i)_{i\in\zz/6}$ where $n_i\in\{0,1,2\}$ and where $(-1)^i(n_i-n_{i+3})$ is independent on~$i$.
\end{mathlist}

{\upshape (b).} The length 6 structure sequences satisfying the condition of {\upshape (\ref{tr10})} are precisely
\[ \begin{array}{c@{\hspace{2em}}c@{\hspace{2em}}c}
 (0, 0, 0, 0, 0, 0) & (0, 0, 2, 0, 0, 2) & (1, 2, 2, 1, 2, 2) \\
 (0, 0, 1, 0, 0, 1) & (0, 2, 0, 2, 0, 2) & (2, 2, 2, 2, 2, 2) \\
 (0, 1, 0, 1, 0, 1) & (0, 2, 2, 0, 2, 2) & (1, 1, 1, 0, 2, 0) \\
 (0, 1, 1, 0, 1, 1) & (1, 1, 2, 1, 1, 2) & (1, 1, 1, 2, 0, 2) \\
 (1, 1, 1, 1, 1, 1) & (1, 2, 1, 2, 1, 2) & (0, 1, 2, 0, 1, 2)
\end{array} \]
up to cyclic permutation and reversing.
\end{theo}

\begin{proof} Define $P(v)$ ($v\in V$) and $\phi_{vs}\col P(v)\ra P(vs)$ uniquely by (\ref{tr4}), (\ref{tr6}), (\ref{tr18}). By the definition of realisations of coloured graphs, $\Gamma$ is a $(2,3,\infty)$-graph if and only if (\ref{tr1}), (\ref{tr5}) and (\ref{tr3}) hold.

Let $P^*(v)$ be the dual to $P(v)$. Let $\<\cdot,\cdot\>\col P(v)\times P^*(v)\ra\rr$ be the natural pairing and let $\{p^*(v,s)\mid s\in S\}$ be the dual basis of $P^*(v)$ defined by
\[ \big\<p(v,s),p^*(v,t)\big\>=\left\{ \begin{array}{@{}ll} 1\qquad &\text{if $s=t$,} \\ 0 & \text{otherwise.} \end{array} \right. \]
Then $\phi_{vs}^{-1}$ induces a map $\phi_{vs}^*\col P^*(v)\ra P^*(w)$. For all $v\in V$ and all distinct $s,t\in S$ we have
\begin{align} p^*(v,s)\, \phi_{vs} &= - p^*(vs,s), \label{tr45}  \\
p^*(v,t)\, \phi_{vs} &= p^*(vs,t) + n(v;s,t)\, p^*(vs,s). \notag \end{align}

Let (\ref{tr44}${}_{k=2}$) and (\ref{tr44}${}_{k=3}$) denote the relevant special cases of the following statement.
\begin{mathlist}

\item \label{tr44} Let $v_0$, $s$, $t$ be such that $m(v_0;s,t)=k$. Let $\{v_i\mid i\in\zz/2k\}$ be the $\{s,t\}$-residue through $v_0$, and $v_is_i=v_{i+1}$ for all $i$, and $s_i=s$ for even $i$ and $s_i=t$ for odd $i$. Then $\phi_{v_1s_1}\cdots\phi_{v_ns_n}=1$.

\end{mathlist}
Then (\ref{tr1}) is equivalent to (\ref{tr44}${}_{k=2}$)  and  (\ref{tr44}${}_{k=3}$). We begin by proving that (\ref{tr44}${}_{k=3}$) is equivalent to (\ref{tr10}) if $\#S\geq3$. Let $s,t,v_i,s_i$ be as in (\ref{tr44}${}_{k=3}$) and let $r\in S\smallsetminus\{s,t\}$. Define the rows of vectors
\begin{align}
f_i:=\big( p^*(v_i,s),p^*(v_i,t),p^*(v_i,r) \big) &&& \text{if $i$ is even,} \label{tr7} \\
f_i:=\big( p^*(v_i,t),p^*(v_i,s),p^*(v_i,r) \big) &&& \text{if $i$ is odd.} \label{tr8}
\end{align}
After interchanging $s,t$ if necessary, we have $f_{i}\, \phi_{v_is_i}=f_i\, M(n_i)$ for all $i$, where $\phi_{v_is_i}$ acts componentwise and
\[ M(n)=\begin{pmatrix} 0 & -1 & 0 \\ 1 & 1 & 0 \\ 0 & n & 1 \end{pmatrix}. \]
We have
\begin{align} & M(n_1)\, M(n_2)\, M(n_3)
= \begin{pmatrix} 0 & -1 & 0 \\ 1 & 1 & 0 \\ 0 & n_1 & 1 \end{pmatrix}
\begin{pmatrix} 0 & -1 & 0 \\ 1 & 1 & 0 \\ 0 & n_2 & 1 \end{pmatrix}
M(n_3) \notag \\[1ex]
& = \begin{pmatrix} -1 & -1 & 0 \\ 1 & 0 & 0 \\ n_1 & n_1+n_2 & 1 \end{pmatrix}
\begin{pmatrix} 0 & -1 & 0 \\ 1 & 1 & 0 \\ 0 & n_3 & 1 \end{pmatrix}
= \begin{pmatrix} -1 & 0 & 0 \\ 0 & -1 & 0 \\ n_1+n_2 & n_2+n_3 & 1 \end{pmatrix} \label{tr11} \end{align}
which is an involution. It follows that\vspace*{1.5ex}
\begin{align*} \smash{\prod_{i=1}^6 M(n_i)=1} &\Longleftrightarrow M(n_1)\, M(n_2)\, M(n_3) =M(n_4)\, M(n_5)\, M(n_6) \\ & \Longleftrightarrow
\text{$(n_4,n_5,n_6)=(n_1,n_2,n_3)+(k,-k,k)$ for some $k$} \\ & \Longleftrightarrow \text{$(-1)^i(n_i-n_{i+3})$ is independent of $i$}. \end{align*}
We have proved that (\ref{tr44}${}_{k=3}$) is equivalent to (\ref{tr10}) if $\#S\geq3$. In case $\#S<3$ the proof is the same as above except that $r$ is absent, that is, the last row and column of $M(n)$ are removed.

Next we prove that (\ref{tr44}${}_{k=2}$) is equivalent to (\ref{tr9}). Let $s,t,v_i,s_i$ be as in (\ref{tr44}${}_{k=2}$) and let $r\in S\smallsetminus \{s,t\}$. As before, define $f_i$ by (\ref{tr7}), (\ref{tr8}). After interchanging $s,t$ if necessary, we have $f_{i}\, \phi_{v_is_i}=f_i\, L(n_i)$ for all $i$ where
\[ L(n)=\begin{pmatrix} 0 & -1 & 0 \\ 1 & 0 & 0 \\ 0 & n & 1 \end{pmatrix}. \]
Now
\[ L(n_1)\, L(n_2)=
\begin{pmatrix} 0 & -1 & 0 \\ 1 & 0 & 0 \\ 0 & n_1 & 1 \end{pmatrix}
\begin{pmatrix} 0 & -1 & 0 \\ 1 & 0 & 0 \\ 0 & n_2 & 1 \end{pmatrix}
=
\begin{pmatrix} -1 & 0 & 0 \\ 0 & -1 & 0 \\ n_1 & n_2 & 1 \end{pmatrix} 
\]
from which it readily follows that
\[ L(n_1)\cdots L(n_4)=1\Longleftrightarrow (n_1,n_2,n_3,n_4)=(n_1,n_2,n_1,n_2). \]
This proves that (\ref{tr44}${}_{k=2}$) is equivalent to (\ref{tr9}). (The case $\#S<3$ is again a consequence of the same computation). Therefore, (\ref{tr1}) is equivalent to (\ref{tr9}) and (\ref{tr10}). Assume now (\ref{tr1}). It remains to prove  (\ref{tr5}) and (\ref{tr3}).

Condition  (\ref{tr5}) states that neighbouring open chambers are disjoint; it holds because of the negative sign in (\ref{tr45}).

Finally, we prove (\ref{tr3}) in the case where $m(R;s,t)=3$, leaving the other cases to the reader. On writing $W:=\sum_{u\in S\smallsetminus\{s,t\}} \rr\, p(v_1,u)$ and $g_i:=( p(v_i,s)+W, p(v_i,t)+W)$ we need to prove $g_i+g_{i+3}=0$. Well, we have $g_{i+3}=g_i K$ where $K$ is the transpose of $M(n_1)\,M(n_2)\,M(n_3)$ (which we computed in (\ref{tr11})) without the last row and last column. This proves the promised case of (\ref{tr3}).

This finishes the proof of (a). Part (b) is straightforward.\qed\end{proof}

\afterblock

The {\em rank\,} of a coloured graph $(V,S,m)$ is defined to be $\#S$. By a {\em $(2,3)$-graph\,} we mean a $(2,3,\infty)$-graph $(V,S,m)$ such that $m(v;s,t)\in\{2,3\}$ for all $v,s,t$. We aim to classify the $(2,3)$-graphs of rank~$3$.

The {\em product\,} of two coloured graphs $(V_i,S_i,m_i)$ ($i=1,2$) is defined to be $(V_1\times V_2,S_1\sqcup S_2,m)$ ($\sqcup$ is disjoint union) where
\begin{align*}
(v_1,v_2)\, s_1 &=(v_1s_1,v_2) && \text{if $v_i\in V_i$ for all $i$ and $s_1\in S_1$},\\
(v_1,v_2)\, s_2 &=(v_1,v_2s_2) && \text{if $v_i\in V_i$  for all $i$ and $s_2\in S_2$},\\
m\big( (v_1,v_2); s,t\big) &= m_i(v_i;s,t) && \text{if $s,t\in S_i$}, \\
m\big( (v_1,v_2); s_1,s_2\big) &= 2 && \text{if $s_i\in S_i$ for all $i$}.
\end{align*}
A coloured graph is {\em irreducible\,} if it is not isomorphic to a product of two coloured graphs of positive rank. It is clear that the product of two $(2,3,\infty)$-graphs is again a $(2,3,\infty)$-graph.

\begin{prop} \label{tr49} Up to isomorphism there are just three irreducible $(2,3)$-graphs of rank $3$:
\[ A_3,\ \tilde{A}_2,\ A(3,7). \]
Here, $A_3$, $\tilde{A}_2$ are the usual names\,%
\footnote{For names of Coxeter groups, see \cite[section~6.4.1]{bou}, \cite[2.4]{hum}.}
of Coxeter groups while $A(3,7)$ is defined\,%
\footnote{In \cite{gru1} this arrangement is called $A_1(7)$ and in \cite{gru2} it is $A(7,1)$.}
in figure~\ref{tr61}(b) (more precisely, it is the $(2,3)$-graph dual to the figure) and which can also be defined by the hyperplane arrangement 
\[ xyz(x+y)(y+z)(z+x)(x+y+z)=0. \]
Moreover, $\tilde{A}_2$ is infinite while the other two are finite.
\end{prop}

\begin{proof} Using theorem~\ref{amcg6} this is an easy exercise involving drawings of graphs, and is left to the reader.\qed
\end{proof}

\section{An example} \label{tr58}

\subsection{An extension of the symmetric group}

From now on we fix an integer $n\geq 0$. Let $G_n$ be the free monoid on a set $T_n\subset G_n$ of $\binom{n+1}{2}$ elements written
\[ T_n=\Big\{ t(a,b)=\binom{a}{b} \ \Big|\ a,b\in\{0,1,\ldots,n\},\ a< b \Big\}. \]

A subset $R\subset G_n$ is {\em closed under cyclic permutations\,} if for all $a,b\in G_n$, if $ab\in R$ then $ba\in R$. (We call $ba$ a {\em cyclic permutation\,} of $ab$).

We define $Q_n\subset G_n$ to be the smallest subset, closed under cyclic permutations, containing
\[ \binom{a}{b}\binom{c}{d}\binom{a}{b}\binom{c}{d} \]
whenever $0\leq a< b\leq c< d\leq n$;
\[ \binom{a}{b}\binom{a+x}{b-y}\binom{a}{b}\binom{a+y}{b-x} \]
whenever $x,y\geq0$ and $0\leq a<a+x+y<b\leq n$; and
\[ \binom{a}{b-z} \binom{a+y}{b} \binom{a}{b-x} \binom{a+z}{b} \binom{a}{b-y} \binom{a+x}{b} \]
whenever $x,y,z> 0$ and $0\leq a\leq a+x+y+z=b\leq n$.

In order to motivate the definition of $Q_n$, note that the action of $G_n$ on $\{1,\ldots,n\}$ defined by
\[ \binom{a}{b}(x)=\left\{ \begin{array}{@{}lll@{}} a+b+1-x &\hspace{2mm}& \text{if $a+1\leq x\leq b$,} \\[.5ex] x && \text{otherwise} \end{array} \right. \]
has the property that the elements of $Q_n$ act trivially.

Let $K_n$ be the group presented by the generating set $T_n$ and relations $s^2=1$ for all $s\in T_n$ and the relations in $Q_n$. One of our aims is to show that $K_n$ is naturally the vertex set of a $(2,3)$-graph.

\subsection{Admissible graphs}

We observe now:
\begin{mathlist} \raggedright

\item \label{tr29} For all distinct $a,b\in T_n$ there are unique $k\in\{2,3\}$ and $c_3,\ldots,c_{2k}\in T_n$ such that $ab(c_3\cdots c_{2k})\in Q_n$. Also, $a^2G_n\cap Q_n=\varnothing$ for all $a\in T_n$.

\item \label{tr30} The set $Q_n$ is also invariant under reversal, that is, under the anti-automorphism of $G_n$ which fixes every element of $T_n$.

\end{mathlist}

\begin{defi} We define an action $T_n\times G_n\ra T_n$ written $(a,b)\mapsto a*b$ as follows. Firstly, $a*a=a$ for all $a\in T_n$. Let $a,b,c\in T_n$ and assume that $Q_n$ meets $abcG_n$. Then $a*b=c$. 

Note that this is well-defined by (\ref{tr29}). Also note that $(a*b)*b=a$ for all $a,b\in T_n$ by (\ref{tr30}).
\end{defi}

\begin{defi} For any set $I$, we define $U_I$ to be the set of injective maps $I\ra T_n$. Recall that $F_I$ is the free monoid on $I$. We define an action $U_I\times F_I\ra U_I$ written $(u,g)\mapsto u \td g$ as follows. Let $u\in U_I$, $s\in I$.
\begin{mathlistn}

\item We put $[u\td s](s)=u(s)$.

\item Let $t\in I\smallsetminus\{s\}$ and put $a:= u(t)$, $b:= [u\td t](s)$. Then $[u\td ts](t)$ is $a*b$, that is, the unique $c\in T_n$ such that $abcG_n\cap Q_n\neq\varnothing$.\qed\end{mathlistn}
\end{defi}

Recall that a groupoid is a category all of whose morphisms are isomorphisms. If $X,Y$ are objects of a category $C$ we write $C(X,Y)$ for the set of morphisms of $C$ from $X$ to $Y$. All our categories are on the right, that is, the composition $C(X,Y)\times C(Y,Z)\ra C(X,Z)$ is written $(f,g)\mapsto fg$ (rather than $gf$).

\begin{defi} For any set $I$, we define the groupoid $R_I$ with object set $U_I$ by the presentation with generators
\[ \begin{pmatrix} u\\s\\u\td s \end{pmatrix} \in R(u,u\td s) \qquad\text{whenever }u\in U_I,\ s\in I \]
and relations
\begin{gather*} \begin{pmatrix} u_0\\s\\u_1 \end{pmatrix}
\begin{pmatrix} u_1\\t\\u_2 \end{pmatrix}
\begin{pmatrix} u_2\\s\\u_3 \end{pmatrix}
\begin{pmatrix} u_3\\t\\u_4 \end{pmatrix}
\cdots
\begin{pmatrix} u_{2k-2}\\s\\u_{2k-1} \end{pmatrix}
\begin{pmatrix} u_{2k-1}\\t\\u_{2k} \end{pmatrix}  \end{gather*}
whenever $u_0=u_{2k}$ and either $s=t$ or
\be \label{tr31} h:=u_0(s)\, u_1(t)\,  u_2(s)\, u_3(t) \cdots  u_{2k-2}(s)\, u_{2k-1}(t) \ee
is a power of an element of $Q_n$.\qed\end{defi}

For every $(u_0,s,t)$ the possible values of $k$ in the above are determined in the following easy result.

\begin{lemm} \label{tr46} Let $u_0\in U_I$, and let $s,t\in I$ be distinct. Define $u_i$ ($i\in\zz$) by $u_{2i-1}\td t=u_{2i}=u_{2i+1}\td s$ for all $i$. Let $g$ be the unique element of $u_0(s)\, u_1(t)\, G_n\cap Q_n$ and let $2p\in\{4,6\}$ be its length. Let $q\in\zz_{\geq0}$ and put $k=pq$. Define $h$ by (\ref{tr31}). Then $h=g^q$ and
\be u_{2k}(r)=u_0(r)*g^q. \label{tr47} \ee
Therefore, $u_0=u_{2k}$ $\Leftrightarrow$ $a*g^q=a$ for all $a\in u_0(I)$.
\end{lemm}

\begin{proof} That $h=g^q$ follows readily from the definitions. For even $i$ we have
\begin{align*}
u_i(r)*u_i(s) &= [u_i\td r](r)*[u_i\td rr](s) \\ & = [u_i\td rrs](r) = [u_i\td s](r) = u_{i+1}(r).
\end{align*}
Likewise, for odd $i$ we have $u_i(r)*u_i(t)=u_{i+1}(r)$. By an obvious induction we find that $u_{2k}(r)=u_0(r)*h$ which proves (\ref{tr47}).\qed\end{proof} 

\begin{defi} Let $u\in U_I$. We define a coloured graph $\Gamma(u)=(V,I,m)$ called an {\em admissible graph\,} as follows. Firstly, $V:=R_I(u,-)$, the set of morphisms in $R_I$ from $u$ to any object. The action $V\times S\ra V$ is defined by
\[ (v,s)\mapsto vs:=v\circ\text{\footnotesize$\begin{pmatrix} u_0\\s\\u_0\td s\end{pmatrix}$} \qquad \text{whenever $v\in R_I(u,u_0)$}. \]
We define $m$ as follows. Use the notation of lemma~\ref{tr46} and let $v\in R_I(u,u_0)$. We define $m(v;s,t)$ to be the least $k>0$ divisible by $p$ such that $u_{2k}=u_0$.\qed\end{defi}

It is clear that $\Gamma(u)$ is a coloured graph. Notice that  $\Gamma(u)$ has a natural base vertex $1_u\in R_I(u,u)$. Note that if $u_1$, $u_2$ are isomorphic objects of $R_I$ then there is an isomorphism of coloured graphs $\Gamma(u_1)\ra\Gamma(u_2)$ (preserving $I$ pointwise) but it may not respect the base points.

\begin{rema} Let $I\subset J$ be sets, $v\in U_J$ and $u:=u|_I\in U_I$. It will follow from later results that there exists a unique injective map of coloured graphs (in the obvious sense) $\Gamma(u)\ra\Gamma(v)$ preserving colours and base points. The injectivity is essentially a consequence of theorem~\ref{6}. However, it is not even clear at this stage that such a map exists at all, because the values of $m$ in $\Gamma(v)$ might be greater than those in $\Gamma(u)$. We shall find this unexpected aspect of admissible graphs helpful in the proof of theorem~\ref{tr55}.\end{rema}

\subsection{Equivalence relations}

Recall that $T_n=\{t(a,b)\mid 0\leq a<b\leq n\}$.

\begin{defi} \label{tr51} (a). For a subset $A\subset\{0,1,\ldots,n\}$ we define $T(A):=\{t(a,b)\mid a,b\in A,\ a<b\}\subset T_n$.

(b). Let $u\in U_I$. The {\em support\,} of $u$ is defined to be $\supp(u):=\{ a,b \mid t(a,b)\in u(I)\}$, that is, the smallest $A$ such that $u(I)\subset T(A)$.

(c). Let $u_1,u_2\in U_I$ and write $A_i=\supp(u_i)$ ($i\in\{1,2\}$). We write $u_1\sim u_2$ if there exists a map $f\col A_1\ra A_2$ which is either an increasing bijection or a decreasing one, and $u_2=g\circ u_1$ where $g\col T(A_1)\ra T(A_2)$ is defined by $g(t(a,b))=t(fa,fb)$.

(d). For the sake of question~\ref{tr60}, we include the following definition. Let $u_1,u_2\in U_I$ and suppose $A=\supp(u_1)=\supp(u_2)$. By a cyclic permutation of $A$ we mean a power of the permutation of $A$ which takes every non-maximal element of $A$ to the next bigger element of $A$. We say that $u_1$ is a {\em cyclic permutation\,} of $u_2$ if there exists a cyclic permutation $f$ of $A$ such that $u_2=g\circ u_1$ where $g\col T(A)\ra T(A)$ is defined by $g(t(a,b))=t(fa,fb)$.\qed\end{defi}

Clearly, $\sim$ is an equivalence relation on $U_I$.

\begin{lemm} \label{tr48} Let $u_1,u_2\in U_I$ be such that $u_1\sim u_2$. 

{\rm (a).} Then $u_1\td g\sim u_2 \td g$ for all $g\in F_I$.

{\rm (b).} Write $E_j:=u_j\td F_I$. Then there is a unique isomorphism of $F_I$-sets $f\col E_1\ra E_2$ (that is, a bijection such that $f(u\td g)=(fu)\td g$ for all $u\in E_1$, $g\in F_I$) such that $f(u_1)=f(u_2)$.

{\rm (c).} For $j\in\{1,2\}$, let $R_j\subset R_I$ be the component of $u_j$, that is, the biggest subcategory of $R_I$ whose object set is $E_j$. Then there is a unique isomorphism of categories $g\col R_1\ra R_2$ such that $g(u)=f(u)$ for all objects $u$ ($f$ as in (b)) and
\[ g
\text{\footnotesize$
\begin{pmatrix} u_3\\s\\u_4 \end{pmatrix}$}=
\text{\footnotesize$
\begin{pmatrix} f(u_3)\\s\\f(u_4) \end{pmatrix}$} \]
whenever the left hand side has a meaning.

{\rm (d).} There is a unique isomorphism $\Gamma(u_1)\ra \Gamma(u_2)$ of pointed coloured graphs which preserves $I$ pointwise.
\end{lemm}

\begin{proof} Easy and left to the reader.\qed\end{proof}

\newcommand{\trav}{\footnotesize}

\newcommand{\traw}{\psdots[dotstyle=o, dotsize=3.5pt](w1)(w2)(w3)(w4)(w5)}

\newcommand{\trbd}{\psset{linewidth=.5pt}} 
\newcommand{\trbe}{\psset{linewidth=1pt}} 

\newcommand{\trax}[1]{%
\psset{unit=8.5mm} \SpecialCoor \degrees[20]
\psframebox[linestyle=none, framesep=0mm]
{\pspicture[.5](-1.2,-1.6)(1.2,1.3)$ 
{\trbd \psarc(0,0){1}{7}{3}}
\pnode(1;7){w1}
\pnode(1;11){w2}
\pnode(1;15){w3}
\pnode(1;19){w4}
\pnode(1;3){w5} \trav #1
$\endpspicture}}

\newcommand{\trbt}{\normalsize \everypsbox{\scriptstyle} \psset{labelsep=4pt}}

\newcommand{\trba}{{%
\trbt
\uput[7]{0}(1;7){0}
\uput[11]{0}(1;11){1}
\uput[15]{0}(1;15){3}
\uput[19]{0}(1;19){6}
\uput[3]{0}(1;3){7}}}

\newcommand{\trbh}{{%
\trbt
\uput[7]{0}(1;7){0}
\uput[11]{0}(1;11){1}
\uput[15]{0}(1;15){2}
\uput[19]{0}(1;19){5}
\uput[3]{0}(1;3){7}}}

\newcommand{\trbc}{{%
\trbt
\uput[-90]{0}(w1){0}
\uput[-90]{0}(w2){1}
\uput[-90]{0}(w3){3}
\uput[-90]{0}(w4){6}
\uput[-90]{0}(w5){7}}}

\newcommand{\trbb}{%
\psframebox[linestyle=none, framesep=0mm]{%
\psset{unit=.35mm, arcangle=90, ncurv=.8} \trbe
\SpecialCoor \trav
\pspicture[.5](0,-30)(100,40)$
{\trbd \psline(0,0)(100,0)}
\pnode(10,0){w1}
\pnode(30,0){w2}
\pnode(50,0){w3}
\pnode(70,0){w4}
\pnode(90,0){w5}
\ncarc{w1}{w3} \ncput[npos=.35, framesep=1pt]{\trcd r} \ncarc{w2}{w4} \ncput[npos=.65, framesep=1pt]{\trcd s} \ncarc[ncurv=1]{w2}{w5} \ncput[npos=.6, framesep=1pt]{\trcd t} \traw \trbc
$\endpspicture}}

\newcommand{\tray}[1]{{\trax{#1} \psset{arcangle=.8, 
nodesep=0pt} \trbe \trav
\ncline{w1}{w3} \ncput[npos=.22, framesep=1pt]{\trcd r} \ncline{w2}{w4} \ncput[npos=.7, framesep=1pt]{\trcd s} \ncline{w5}{w2} \ncput[npos=.3, framesep=1pt]{\trcd t} \traw}}

\newcommand{\trbi}[1]{{\trax{#1} \psset{arcangle=.8, 
nodesep=0pt} \trbe \trav
\ncline{w4}{w1} \ncput[npos=.78, framesep=1pt]{\trcd r} \ncline{w3}{w5} \ncput[npos=.3, framesep=1pt]{\trcd s} \ncline{w5}{w2} \ncput[npos=.7, framesep=1pt]{\trcd t} \traw}}

\newcommand{\traz}[1]{{\trax{#1} \psset{arcangle=.8, 
nodesep=0pt} \trbe \trav
\ncline{w1}{w3} \ncline{w2}{w4} \ncline{w5}{w2} \traw}}

\newcommand{\trcd}[1]{\pscirclebox*
[framesep=1pt]{#1}}                   

\afterblock
Let $\approx$ be the equivalence relation on $U_I$ generated by $\sim$ defined in definition~\ref{tr51} and $\cong$ (isomorphism in the groupoid $R_I$). 

Let $\approx_s$ be the equivalence relation on $U_I$ generated by $\approx$ and the graph of the symmetric group on $I$. In other words, $u_1\approx_s u_2$ if and only if $u_1\approx u_2\circ s$ for some permutation $s$ of $I$. Define $\sim_s$ and $\cong_s$ likewise.

\begin{figure}[h]
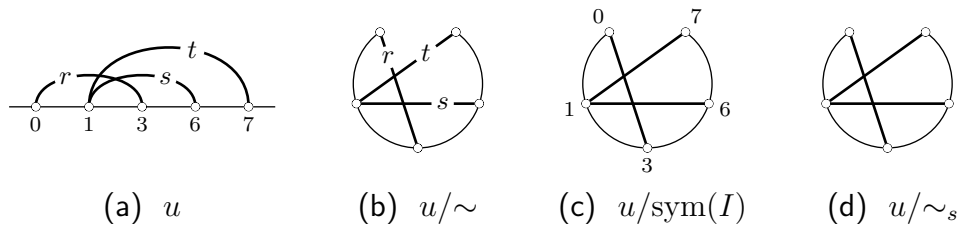

\begin{gather*} \begin{tabular}{c@{\hspace{8mm}}c@{\hspace{8mm}}c@{\hspace{8mm}}c} \trbb & $\tray{}$ & \traz{\trba} & \traz{}   \\ (a) \ $u$ & (b) \  $u/{\sim}$ & (c)  \ $u/{\sym(I)}$ & (d) \  $u/{\sim_s}$  \end{tabular} \end{gather*}
\caption{\label{tr63} Vertices of admissible graphs. Let $I=\{r,s,t\}$ have $3$ elements. Part (a) shows a picture of the object $u\in U_I$ defined by 
$u(r)=t(0,3)$, $u(s)=t(1,6)$, $u(t)=t(1,7)$. The picture in (a) is flat but we usually prefer the (equivalent) curled up version of (b)--(d). In (b) we see the $\sim$-class of $u$. The precise values $0,1,3,6,7$ are forgotten but their ordering is not as it is still shown in the picture. In (c) we divide out the symmetric group on $I$ and in (d) we divide out $\sim_s$.}
\end{figure}

\begin{figure}[h]
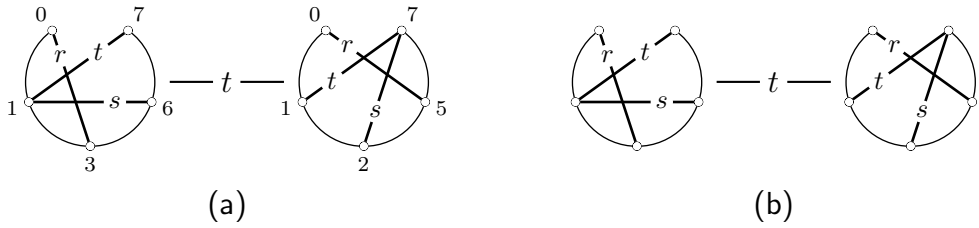

\begin{gather*} \begin{tabular}{@{}c@{\hspace{15mm}}c@{}}
\psmatrix[colsep=15mm] $\tray{\trba}$ & $\trbi{\trbh}$ \endpsmatrix \ncline{1,1}{1,2} \ncput*{$t$}
&
\psmatrix[colsep=15mm] $\tray{}$ & $\trbi{}$ \endpsmatrix \ncline{1,1}{1,2} \ncput*{$t$} \\ (a) & (b)
\end{tabular} \end{gather*}
\caption{\label{tr64} Edges of admissible graphs. Suppose that the left hand side in (a) depicts some $u\in U_I$ with $I=\{r,s,t\}$. Then the right hand side is $u\td t$. Part (b) is obtained from (a) by taking $\sim$-classes --- we know by lemma~\ref{tr48} that coloured graphs survive division by~$\sim$.}
\end{figure}

It is natural to draw pictures of objects of $R_I$. The convention is easily understood from figure~\ref{tr63} which shows pictures for an element of $U_I$ and some of its equivalence classes, and figure~\ref{tr64} which shows pictures for edges in admissible graphs.

\subsection{Main result}

\begin{exam} \label{tr52} We now have a detailed look at three rank~$3$ admissible graphs. Our understanding of them will be crucial in the case-by-case proof of theorem~\ref{tr55}.

(a). One rank~$3$ admissible graph $\Gamma$ is given in figure~\ref{tr39}. You should verify it. Note that the vertices can be taken to be $\sim$-classes by lemma~\ref{tr48}. The verification is helped by the schematic version of the graph in figure~\ref{tr40} and the order $6$ automorphism group (which fixes the $2$-residue numbered $1$).

An observation which will be important in the proof of theorem~\ref{tr55} is that $\Gamma$ is a $(2,3)$-graph.  Indeed, it is isomorphic to the Coxeter graph of type~$A_3$.

(b). Figure~\ref{tr61}(a) shows part of another rank~$3$ admissible graph $\Gamma$. Convince yourself that it is correct. The dashed triangle is precisely $1/8$ of the whole graph. The automorphism group of $\Gamma$ is of order $8$ and generated by the reflections in the edges of the dashed triangle. 

Again, we observe that $\Gamma(u)$ is a $(2,3)$-graph (use theorem~\ref{amcg6} or proposition~\ref{tr49}). In the classification of rank~$3$ $(2,3)$-graphs (proposition~\ref{tr49}) we said that it is of type $A(3,7)$. As every $(2,3)$-graph, it has a realisation as a hyperplane arrangement. This arrangement is shown in figure~\ref{tr61}(b), which also serves to give a full picture rather than $1/8$ of it.

(c). Let $u\in U_I$ be such that $u(I)=\{ t(0,2), t(1,3), t(2,4) \}$. Then $u$ is a single isomorphism class in $R_I$ and one easily deduces that $\Gamma(u)$ must be a Coxeter coloured graph. Indeed it is of type $A_3$ and again it is a $(2,3)$-graph.
\end{exam}

\begin{defi} Let $u\in U_I$. We call $u$ {\em reducible\,} if $I$ can be written as the union of two non-empty disjoint sets $A,B$ such that for all $(a,b)\in A\times B$ there exist $x,y\in T_n$ such that $u(a)\, u(b)\, x\, y\in Q_n$. Otherwise it is called {\em irreducible}.
\end{defi}

\begin{lemm} \label{tr50} Let $u\in U_I$. If $u$ is reducible then $\Gamma(u)$ is reducible as a coloured graph.
\end{lemm}

\begin{proof} Left to the reader.\qed\end{proof}


\newcommand{\trbk}{\psdots[dotstyle=o, dotsize=3.5pt](w1)(w2)(w3)(w4)(w5)}

\newcommand{\trbj}{%
\psset{unit=6mm} \SpecialCoor \degrees[20] \pspicture(-1,-1.2)(1,1) \psarc[linewidth=.5pt](0,0){1}{7}{3}
\pnode(1;7){w1}
\pnode(1;11){w2}
\pnode(1;15){w3}
\pnode(1;19){w4}
\pnode(1;3){w5}
\endpspicture}

\newcommand{\trbl}{{\trbj 
\ncline{w1}{w3} \ncline{w2}{w4} \ncline{w5}{w3} \trbk}}

\newcommand{\trbm}{{\trbj 
\ncline{w1}{w3} \ncline{w2}{w4} \ncline{w5}{w2} \trbk}}

\newcommand{\trbn}{{\trbj 
\ncline{w1}{w3} \ncline{w2}{w5} \ncline{w4}{w1} \trbk}}


\newcommand{\trbo}{\psdots[dotstyle=o, dotsize=3.5pt](w1)(w2)(w3)(w4)(w5)(w6)}

\newcommand{\trbp}{%
\psset{unit=6mm} \SpecialCoor \degrees[6] \pspicture(-1,-1.2)(1,1) \psarc[linewidth=.5pt](0,0){1}{2}{1}
\pnode(1;2){w1}
\pnode(1;3){w2}
\pnode(1;4){w3}
\pnode(1;5){w4}
\pnode(1;0){w5}
\pnode(1;1){w6}
\endpspicture}

\newcommand{\trbq}{{\trbp
\ncline{w1}{w4} \ncline{w2}{w5} \ncline{w6}{w3} \trbo}}

\newcommand{\trbr}{{\trbp
\ncline{w1}{w3} \ncline{w2}{w5} \ncline{w6}{w4} \trbo}}

\newcommand{\trbs}{{\trbp
\ncline{w1}{w4} \ncline{w2}{w6} \ncline{w5}{w3} \trbo}}



\newcommand{\trab}{\psdots[dotstyle=o, dotsize=3.5pt](w1)(w2)(w3)(w4)(w5)}

\newcommand{\traa}{%
\psset{unit=8mm} \SpecialCoor \degrees[20] \pspicture(0,0)(0,0)$ \psarc[linewidth=.5pt](0,0){1}{7}{3}
\pnode(1;7){w1}
\pnode(1;11){w2}
\pnode(1;15){w3}
\pnode(1;19){w4}
\pnode(1;3){w5}
$\endpspicture}

\newcommand{\trbw}{\trbu \traa \psset{arcangle=.8, 
nodesep=0pt, linewidth=1pt}}

\newcommand{\trbu}{\small}       
\newcommand{\trae}[1]{\pscirclebox*
[framesep=.5pt]{\scriptstyle #1}}       
\newcommand{\trbv}{\small}              

\newcommand{\trac}[2]{{\trbw
\ncarc{w1}{w3} \ncput[npos=.23, framesep=1pt]{\trae {s}} \ncarc{w2}{w4} \ncput[npos=.7, framesep=1pt]{\trae {#1}} \ncarc{w5}{w2} \ncput[npos=.3, framesep=1pt]{\trae {#2}} \trab}}

\newcommand{\trad}[2]{{\trbw
\ncarc{w4}{w1} \ncput[npos=.77, framesep=1pt]{\trae {s}} \ncarc{w3}{w5} \ncput[npos=.3, framesep=1.5pt]{\trae {#1}} \ncarc{w5}{w2} \ncput[npos=.7, framesep=1pt]{\trae {#2}} \trab}}

\newcommand{\traf}{\trac{r}{t}}
\newcommand{\trag}{\trac{t}{r}}
\newcommand{\trah}{\trad{r}{t}}
\newcommand{\trai}{\trad{t}{r}}

\newcommand{\traj}[6]{#5%
\SpecialCoor \pspicture(-14,-25)(14,25)$
\rput(0,-8){1} 
\pnode(4,4){v1}    \pnode(-4,4){v2}
\pnode(-4,-4){v3}  \pnode(4,-4){v4}
\pnode(4,12){v5}  \pnode(-4,12){v6}
\pnode(-4,-12){v7}\pnode(4,-12){v8}
\pnode(12,4){v9}  \pnode(-12,4){v10}
\pnode(-12,-4){v11}  \pnode(12,-4){v12}
\pnode(4,20){v13}    \pnode(-4,20){v14}
\pnode(-4,-20){v15}  \pnode(4,-20){v16}
\pnode(12,12){v17}    \pnode(-12,12){v18}
\pnode(-12,-12){v19}  \pnode(12,-12){v20}
\pnode(10,25){v21}    \pnode(-10,25){v22}
\pnode(-10,-25){v23}  \pnode(10,-25){v24}
\ncline{v1}{v2} \ncput{\trbx{t}}  
\ncline{v2}{v3} \ncput{\trbx{r}}
\ncline{v3}{v4} \ncput{\trbx{t}}  
\ncline{v4}{v1} \ncput{\trbx{r}}
\ncline{v1}{v5} \ncput{\trbx{s}}  \ncline{v2}{v6} \ncput{\trbx{s}}
\ncline{v3}{v7} \ncput{\trbx{s}}  \ncline{v4}{v8} \ncput{\trbx{s}}
\ncline{v5}{v9} \ncput{\trbx{r}}  \ncline{v6}{v10} \ncput{\trbx{r}}
\ncline{v7}{v11} \ncput{\trbx{r}} \ncline{v8}{v12} \ncput{\trbx{r}}
\ncline{v5}{v13} \ncput{\trbx{t}}
\ncline{v6}{v14} \ncput{\trbx{t}}
\ncline{v7}{v15} \ncput{\trbx{t}}
\ncline{v8}{v16} \ncput{\trbx{t}}
\ncline{v9}{v17} \ncput{\trbx{t}}
\ncline{v10}{v18} \ncput{\trbx{t}}
\ncline{v11}{v19} \ncput{\trbx{t}}
\ncline{v12}{v20} \ncput{\trbx{t}}
\ncline{v13}{v17} \ncput{\trbx{r}}
\ncline{v14}{v18} \ncput{\trbx{r}}
\ncline{v15}{v19} \ncput{\trbx{r}}
\ncline{v16}{v20} \ncput{\trbx{r}}
\ncline{v13}{v14} \ncput{\trbx{s}}
\ncline{v10}{v11} \ncput{\trbx{s}}
\ncline{v15}{v16} \ncput{\trbx{s}}
\ncline{v12}{v9}  \ncput{\trbx{s}}
\ncline{v17}{v21}  \ncput{\trbx{s}}
\ncline{v18}{v22}  \ncput{\trbx{s}}
\ncline{v19}{v23}  \ncput{\trbx{s}}
\ncline{v20}{v24}  \ncput{\trbx{s}}
\ncline{v21}{v22}  \ncput{\trbx{r}}
\ncline{v23}{v24}  \ncput{\trbx{r}}
{#6 \ncdiag[angle=0, linearc=1]{v21}{v24} \ncput{\trbx{t}}
\ncdiag[angle=180, linearc=1]{v22}{v23} \ncput{\trbx{t}}}
\rput{0}(v1){#4} \rput{0}(v2){#4} \rput{0}(v3){#2} \rput{0}(v4){#2} \rput{0}(v5){#3} \rput{0}(v6){#3} \rput{0}(v7){#2} \rput{0}(v8){#2} \rput{0}(v9){#3} \rput{0}(v10){#3} \rput{0}(v11){#4} \rput{0}(v12){#4} \rput{0}(v13){#1} \rput{0}(v14){#1} \rput{0}(v15){#2} \rput{0}(v16){#2} \rput{0}(v17){#1} \rput{0}(v18){#1} \rput{0}(v19){#4} \rput{0}(v20){#4} \rput{0}(v21){#1} \rput{0}(v22){#1} \rput{0}(v23){#3} \rput{0}(v24){#3}
$\endpspicture}
\newcommand{\trat}[1]{\pscirclebox[linewidth=.5pt, fillstyle=solid, fillcolor=white, framesep=1pt]{\scriptstyle #1}}

\newcommand{\trbx}[1]{\pscirclebox*
[framesep=1.5pt]{\scriptstyle #1}}  

\begin{figure}
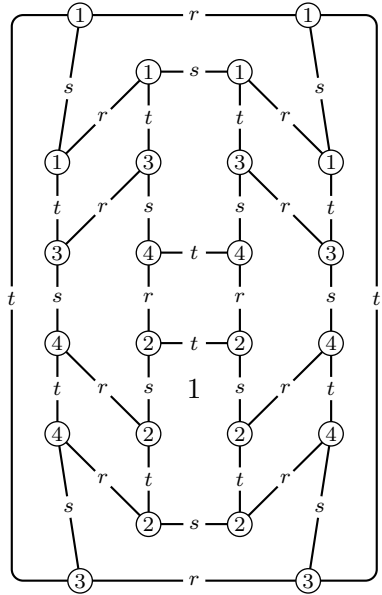

\caption{Schematic version of $\Gamma(u)$ for $u\in L_2\cup L_3$. See figure~\ref{tr39} for a full picture. The automorphism group of this graph has order $6$ and preserves the $2$-residue labelled~$1$.\label{tr40}}
\begin{align*} \small
\psframebox[linestyle=none, framesep=2mm]{\traj 
{\trat{1}} {\trat{2}} {\trat{3}} {\trat{4}}
{\psset{unit=1.5mm, framesep=1.5pt, nodesep=0mm}} {\psset{arm=6}} }
\end{align*}
\end{figure}

\renewcommand{\trbx}[1]{\pscirclebox*
[framesep=2pt]{#1}}                   

\begin{figure}
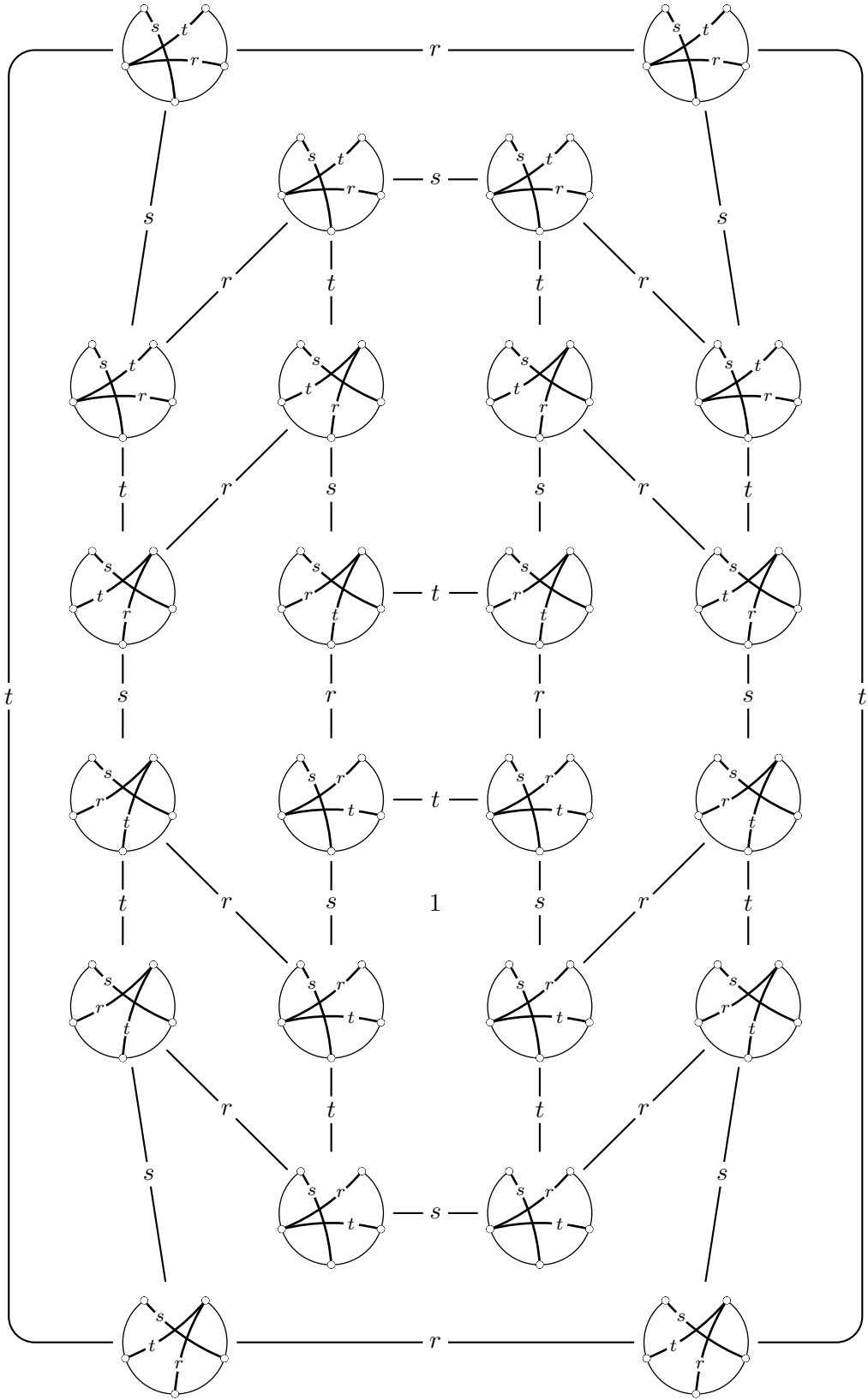
 \begin{center} \caption{The coloured graph $\Gamma(u)$ for $u\in L_2\cup L_3$. A schematic version of it is shown in figure~\ref{tr40}.\label{tr39}} \trbv
\[ \psframebox[linestyle=none, framesep=8mm]{\traj {\traf} {\trag} {\trah} {\trai} {\psset{unit=4.0mm, nodesep=9.5mm}} {\psset{arm=4}} } \] \end{center} \end{figure}


\newcommand{\tram}{\small
\psset{npos=.28, framesep=1.5pt}}       
\newcommand{\trbz}[1]{\pscirclebox*
[framesep=.5pt]{\scriptstyle #1}}       
\newcommand{\trca}{\small}              
\newcommand{\trcb}[1]{\pscirclebox*
[framesep=1pt]{#1}}                     

\newcommand{\trak}{\psdots[dotstyle=o, dotsize=3.5pt](w1)(w2)(w3)(w4)(w5)(w6)}

\newcommand{\tral}{%
\psset{unit=7mm} \SpecialCoor \degrees[12] \pspicture(0,0)(0,0)$ \psarc[linewidth=.5pt](0,0){1}{4}{2}
\pnode(1;4){w1} \pnode(1;6){w2} \pnode(1;8){w3} \pnode(1;10){w4} \pnode(1;12){w5} \pnode(1;14){w6}
$\endpspicture}

\newcommand{\tran}{{\tral \psset{nodesep=0pt, linewidth=1pt} 
{\tram 
\ncline{w1}{w4} \ncput{\trbz{r}}
\ncline{w5}{w2} \ncput{\trbz{s}} 
\ncline{w3}{w6} \ncput{\trbz{t}}}
\trak}}

\newcommand{\trao}{{\tral \psset{nodesep=0pt, linewidth=1pt} {\tram 
\ncline{w1}{w4} \ncput[npos=.5]{\trbz{r}}
\ncline{w5}{w3} \ncput{\trbz{s}} 
\ncline{w2}{w6} \ncput{\trbz{t}}}
\trak}}

\newcommand{\traq}{{\tral \psset{nodesep=0pt, linewidth=1pt} {\tram 
\ncline{w2}{w5} \ncput[npos=.5]{\trbz{s}}
\ncline{w1}{w3} \ncput{\trbz{r}} 
\ncline{w4}{w6} \ncput{\trbz{t}}}
\trak}}

\newcommand{\trar}{{\tral \psset{nodesep=0pt, linewidth=1pt} {\tram 
\ncline{w3}{w6} \ncput[npos=.5]{\trbz{t}}
\ncline{w5}{w1} \ncput{\trbz{r}} 
\ncline{w2}{w4} \ncput{\trbz{s}}}
\trak}}

\newcommand{\tras}[4]{%
\SpecialCoor \psset{unit=2.2pt, nodesep=8mm, framesep=1.5pt}
\psframebox[linestyle=none, framesep=1mm]
{\pspicture(-80,-91)(80,72)$
\rput(-78,70){\text{\normalsize (a)}}
\pnode(0,0){v4}
\pnode(58.3;150){v1} \pnode(58.3;30){v2} \pnode(58.3;-90){v3}
\pnode(34;150){v5} \pnode(34;30){v6} \pnode(34;-90){v7}
\pnode(-65;150){v8} \pnode(-65;30){v9} \pnode(-65;-90){v10}
\pnode(60;125){v11} \pnode(60;5){v12} \pnode(60;245){v13}
\pnode(60;55){v14} \pnode(60;-65){v15} \pnode(60;175){v16}
\pnode(83;150){v17} \pnode(83;30){v18} \pnode(83;-90){v19}
\rput{0}(v4){#1} 
\rput{0}(v5){#4} \rput{0}(v6){#2} \rput{0}(v7){#3}
\rput{0}(v8){#1} \rput{0}(v9){#1} \rput{0}(v10){#1}
\rput{0}(v11){#4} \rput{0}(v12){#2} \rput{0}(v13){#3}
\rput{0}(v14){#2} \rput{0}(v15){#3} \rput{0}(v16){#4}
\rput{0}(v17){#4} \rput{0}(v18){#2} \rput{0}(v19){#3}
{\psset{arcangle=-15, nodesep=0, arrows=c-c} \dashed \ncarc{v2}{v1} \ncarc{v1}{v3} \ncarc{v3}{v2}}
\ncline{v4}{v5} \ncput{\trcb{t}} \ncline{v4}{v6} \ncput{\trcb{r}} \ncline{v4}{v7} \ncput{\trcb{s}}
\ncline[border=3mm]{v5}{v11} \ncput{\trcb{r}}
\ncline[border=3mm]{v6}{v12} \ncput{\trcb{s}}
\ncline[border=3mm]{v7}{v13} \ncput{\trcb{t}}
\ncline[border=3mm]{v6}{v14} \ncput{\trcb{t}}
\ncline[border=3mm]{v7}{v15} \ncput{\trcb{r}}
\ncline[border=3mm]{v5}{v16} \ncput{\trcb{s}}
\ncline{v17}{v11} \ncput{\trcb{s}} \ncline{v18}{v12} \ncput{\trcb{t}} \ncline{v19}{v13} \ncput{\trcb{r}}
\ncline{v18}{v14} \ncput{\trcb{s}} \ncline{v19}{v15} \ncput{\trcb{t}} \ncline{v17}{v16} \ncput{\trcb{r}}
\ncline{v10}{v11} \ncput{\trcb{t}} \ncline{v8}{v12} \ncput{\trcb{r}}  \ncline{v9}{v13} \ncput{\trcb{s}}
\ncline{v10}{v14} \ncput{\trcb{r}} \ncline{v8}{v15} \ncput{\trcb{s}}  \ncline{v9}{v16} \ncput{\trcb{t}}
$\endpspicture}}

\begin{figure} 
\caption{{\bfseries Picture (a).} Part of the admissible graph $\Gamma(u)$ for $u\in L_4\cup L_5\cup L_6$. The dashed triangle is exactly one eighth of it and corresponds to the gray region of (b) and (c). {\bfseries Pictures (b) and (c).} The line arrangement defined by $xyz(x+y)(y+z)(z+x)(x+y+z)=0$. It is dual to the $(2,3)$-graph $A(3,7)$. The gray region is the dashed triangle of 
(a).\label{tr61}}
\begin{gather*} \footnotesize \tras{\tran}{\trao}{\traq}{\trar} 
\\ 
\psframebox[linestyle=none, framesep=0pt]{\begin{array}{llrl}
\psset{unit=8mm, dimen=middle}
\psframebox[linestyle=none, framesep=0mm] 
{\pspicture[.5](-3.5,-3.5)(3.5,3.5)
\pswedge[linestyle=none, fillstyle=solid, fillcolor=lightgray](0,0){1.4142}{0}{90}
\pscircle(0,0){1.4142}
\pscircle(1,1){2} \pscircle(-1,1){2} \pscircle(1,-1){2} \pscircle(-1,-1){2}
\psline(-3.5,0)(3.5,0) \psline(0,-3.5)(0,3.5)
\endpspicture} & \hspace{5mm} &&
\psset{unit=10mm}
\psframebox[linestyle=none, framesep=0mm] 
{\pspicture[.5](-2,-2.5)(2.5,2.5)$
\pswedge[linestyle=none, fillstyle=solid, fillcolor=lightgray](0,0){2.236}{0}{90}
\psline(1,2)(-2,-1) \psline(-1,-2)(2,1)
\psline(1,-2)(-2,1) \psline(-1,2)(2,-1)
\psline(0,-2.236)(0,2.236) \psline(-2.236,0)(2.236,0)
\rput(1.9,1.9){\infty}
$\endpspicture} 
\\ \\ \text{(b) \ Spherical picture of $A(3,7)$.} 
&& \text{(c)} & \text{Projective picture of $A(3,7)$.} \\ &&& \text{The line at infinity is included.}
\end{array}} 
\end{gather*}
\end{figure}

\begin{lemm} \label{tr54} Suppose that $n\geq 4$ and $I$ is a set of $3$ elements.

{\upshape (a).} There are precisely six $\sim_s$-classes $L_1$, \ldots, $L_6$ of irreducible elements in $U_I$. They are given by the following representatives.
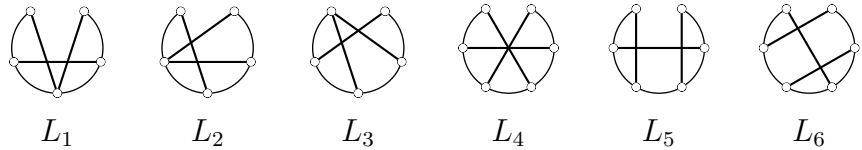
\begin{align*} \psset{linewidth=1pt} && &&
\begin{array}{c} \trbl\\ L_1\end{array} &&
\begin{array}{c} \trbm\\ L_2\end{array} &&
\begin{array}{c} \trbn\\ L_3\end{array} &&
\begin{array}{c} \trbq\\ L_4\end{array} &&
\begin{array}{c} \trbr\\ L_5\end{array} &&
\begin{array}{c} \trbs\\ L_6\end{array} && &&
\end{align*}

{\upshape (b).} The $\approx_s$-classes of irreducible elements in $U_I$ are $L_1$, $L_2\cup L_3$ and $L_4\cup L_5\cup L_6$.
 
{\upshape (c).} Every rank~$3$ admissible graph is a $(2,3)$-graph.
\end{lemm}

\begin{proof} It is easy and left to the reader to prove (a) using lemma~\ref{tr50}. 

Proof of (b). The (connected) graph of example~\ref{tr52}(a) and figure~\ref{tr39} involves $L_2$ and $L_3$ but no others (recall that reflection through a vertical line fixes every $\sim$-class by definition). Therefore $L_2\cup L_3$ is a single $\approx_s$-class. Likewise, the graph of example~\ref{tr52}(b) and figure~\ref{tr61}(a) involves $L_4$, $L_5$ and $L_6$ but no others so $L_4\cup L_5\cup L_6$ is a $\sim_s$-class. Only one $\sim_s$-class $L_1$ remains which must therefore be a $\approx_s$-class as well; we looked at the related admissible graph in example~\ref{tr52}(c).

Proof of (c). By (b) and lemma~\ref{tr50} we know all {\em irreducible\,} rank~$3$ admissible graphs. As we already observed in example~\ref{tr52}, all of them are $(2,3)$-graphs. It is easy and left to the reader to handle the reducible ones.\qed\end{proof}


\begin{theo} \label{tr55} Every admissible graph is a $(2,3)$-graph.
\end{theo}

\begin{proof} Consider an admissible graph $\Gamma(u)=(V,I,m)$, $u\in U_I$. 

First we prove that $m(v;s,t)\in\{2,3\}$ for all $v,s,t$. In lemma~\ref{tr54} we observed this to be true in the rank~$3$ case. By lemma~\ref{tr46}, this implies that $a*g=a$ for all $(a,g)\in T_n\times Q_n$. Using lemma~\ref{tr46} backwards we find that $m(v;s,t)\in\{2,3\}$ for all $v,s,t$.

Recall that a $(2,3)$-graph is just a $(2,3,\infty)$-graph for which $m(v;s,t)$ is never infinite. By theorem~\ref{amcg6} it remains to prove that all structure sequences of $\Gamma(u)$ satisfy (\ref{tr9}) and (\ref{tr10}). But all structure sequences of all admissible graphs occur in rank~$3$ admissible graphs. In lemma~\ref{tr54} we already observed the latter to be $(2,3)$-graphs, in particular, to satisfy the required conditions  (\ref{tr9}) and (\ref{tr10}).\qed\end{proof}

\begin{coro} There exists a faithful linear representation of $K_n$.
\end{coro}

\begin{proof} Let $u\in U_I$ be such that $u\col I\ra T_n$ is surjective. By theorem~\ref{tr55}, $\Gamma(u)$ is a $(2,3)$-graph. It is clear that $K_n$ acts on $\Gamma(u)$. By the unicity of the standard realisation of $(2,3)$-graphs, this action passes to a $K_n$-action on $Q$. The action on $Q$ is faithful because the action on $\Gamma(u)$ is.\qed\end{proof}

\begin{ques} \label{tr60} Recall that in definition~\ref{tr51}(d) we defined cyclic permutations of elements of $U_I$. Observe now that every $u_1\in L_1$ is a cyclic permutation of some $u_2\in L_2$ (see lemma~\ref{tr54}(a) for the classification of rank~$3$ admissible graphs). Also, $\Gamma(u_1)$ and $\Gamma(u_2)$ are isomorphic as coloured graphs because both are of Coxeter type $A_3$ as we saw in example~\ref{tr52}(a) and (c). I don't know if this is a coincidence. Is it true in general that $\Gamma(u_3)$ and $\Gamma(u_4)$ are isomorphic whenever $u_3$ is a cyclic permutation of $u_4$?
\end{ques}

We finish with a result without proof.

\begin{prop} \label{tr59} There are precisely four isomorphism classes of rank~$4$ irreducible finite $(2,3)$-graphs. They are the Coxeter ones $A_4$, $D_4$ and two more named $A(4,13)$, $A(4,15)$. Among them, $D_4$ is the only non-admissible one. Possible choices of $u_{13}$, $u_{15}\in U_I$ such that $A(4,13)=\Gamma(u_{13})$, $A(4,15)=\Gamma(u_{15})$ are as follows.
\begin{align*}
u_{13}(I) &=\{t(0,2),t(0,3),t(0,4),t(1,5)\}, \\
u_{15}(I) &=\{t(0,2),t(0,4),t(1,5),t(3,6)\}.
\end{align*}
Here are possible equations for $A(4,13)$, $A(4,15)$.
\begin{align*}
\text{\upshape A(4,13)\quad} & xyzw(x+y)(y+z)(z+w)(w+y)(y+z+w) \\*
& (x+y+z)(x+y+w)(x+y+z+w)(x+2y+z+w) \\[1ex]
\text{\upshape A(4,15)\quad} & xyzw(x+y)(y+z)(z+w) (w+y)  \\* & (x+y+z) (x+y+w) (y+z+w) (x+2y+z) \\* & (x+y+z+w) (x+2y+z+w) (x+2y+2z+w)
\end{align*}
The Poincar\'e polynomials\,%
\footnote{See \cite[section~2.3]{orl} for the definition of the Poincar\'e polynomial. We use the notation $[n]=1+nt$.}
of $A(3,7)$, $A(4,13)$, $A(4,15)$ are, respectively, $[1][3]^2$, $[1][3][4][5]$, $[1][4][5]^2$.\qed\end{prop}

Daan Krammer, University of Warwick, Mathematics Department, Coventry CV4~7AL, UK, {\tt D.Krammer@warwick.ac.uk} .

\end{document}